\documentclass[a4paper,12pt,oneside]{article}
\usepackage[english]{babel}
\usepackage{amsmath,amsfonts,latexsym,amssymb,amsthm}

\begin{document}
\begin{center}{\Large\bf Large strings of consecutive smooth integers}\\

\vspace{1cm}
Filip Najman
\end{center}

\textbf{Abstract.} In this note we improve an algorithm from a recent paper by Bauer and Bennett for computing a function of Erd\"os that measures the minimal gap size $f(k)$ in the sequence of integers
at  least  one  of  whose  prime  factors  exceeds  $k$. This allows us to compute values of $f(k)$ for larger $k$ and obtain new values of $f(k)$.\\

\par\noindent \textbf{Keywords} Pell equation, Compact representations, Smooth numbers
\par\noindent \textbf{Mathematics subject classification (2000)} 11N25, 11D09.\\

\section{Introduction}

For any integer $m$ we let $P(m)$ be the largest prime factor of $m$ with the convention $P(0)=P(\pm 1)=1$. Let $\Pi_{n,k}$ be the product of $k$ consecutive integers, starting with $n$, i.e.
$$\Pi_{n,k} = n(n+1)\cdots(n+k-1).$$By a theorem of Sylvester (see \cite{syl}) $\Pi_{n,k}$ is divisible by a prime $p>k$ whenever $n>k$, and thus, following Erd\"os \cite{erd}, we define $f(k)$ to be the least integer with the property that
$$P\left(\Pi_{n,f(k)}\right)>k,$$
if $n>k$.
Standard heuristics for the size of gaps between consecutive primes lead one to expect that order of magnitude of $f(k)$ is $(\log k)^2$.

The following table gives known values of $f(k)$:
\begin{center}
\begin{tabular}{|cc|cc|cc|cc|}
\hline
$k$ & $f(k)$ & $k$ & $f(k)$ & $k$ & $f(k)$\\
\hline
$1$ & $1$ & $13-40$ & $6$ & $61-113$& $14$ \\
$2$ & $2$ &  $41-46$ & $7$ & $114$ & $13$ \\
$3-4$ & $3$ & $47-58$ & $8$ & $115-150$ & $12$\\
$5-12$ & $4$ & $59-60$ & $9$ &$151-178$ & $14$\\
\hline
\end{tabular}
\end{center}
The values of $f(k)$ for $k\leq 10$ were computed by Utz \cite{utz} and extended to $k\leq 42$ by Lehmer \cite{leh2}, to $k\leq 46$ by Ecklund and Eggleton \cite{ee}, to $k\leq 73$ by Ecklund, Eggleton and Selfridge \cite{ees1}, \cite{ees2} and finally to $k\leq 178$ by Bauer and Bennett \cite{bb}. Bauer and Bennett in the same paper also disproved an assertion of Utz that $f$ is monotone.

In this paper we compute the values of $f(k)$ for $k\leq 268$. Our results can be summarized in the following theorem.
\newtheorem{tm}{Theorem}
\begin{tm}
For $179 \leq k \leq 268$ the values of $f(k)$ are as follows:
\begin{center}
\begin{tabular}{|cc|}
\hline
$k$ & $f(k)$\\
\hline
$179-222$ & $14$\\
$223-268$ & $16$\\
\hline
\end{tabular}
\end{center}
\end{tm}

\medskip

Note that Bauer and Bennett, although greatly extending the set of $k$ such that $f(k)$ is known, by $105$ values of $k$, did not find any values of $f(k)$ such that $f(k)>14$. In Theorem 1 we find the new largest proven value of $f(k)$, the first after nearly 40 years (the previous being \cite{ees2}).
\section{The algorithm}
Lemher \cite{leh} searched for two consecutive smooth integers, satisfying $P(z(z+1))\leq p_t$, where $p_t$ is the $t$-th prime. One can write $x=2z+1$ and see that $P(z(z+1))\leq p_t$ iff $P(x^2-1)\leq p_t$. Writing $x^2-1=dy^2$, where $d$ is squarefree, leads to the Pell equation \begin{equation}x^2-dy^2=1,\label{pell}
\end{equation}
where $P(dy)\leq p_t$. This implies that the solutions $x+y\sqrt d$ we are searching for are powers of the fundamental solution $u+v\sqrt d$, i.e
$$x+y\sqrt d=(u+v\sqrt d)^n.$$
By classic results on primitive divisors (see \cite{car}) and the work of Lehmer (see \cite{leh3}), if we want $P(y)\leq p_t$, then
$$n\leq \max\left\{  \frac{p_t+1}2, 12\right\}.$$
Thus one needs to consider only finitely many Pell equations and for each one only some of the first solutions.

After finding all pairs of smooth consecutive integers in this way, one can search through the results and find larger strings of consecutive integers.

Bauer and Bennett \cite{bb} improve on this strategy by the following clever argument: for fixed $m\leq 3$ and $t$ suppose we are searching for all integers $n$ satisfying
\begin{equation}
\label{eq1}
P(\Pi_{n,m})\leq p_t.
\end{equation}

We can split the indices $0,\ldots, m-1$ into $\lfloor \frac m 4 \rfloor+ \lfloor \frac {m+1} 4 \rfloor$ disjoint pairs $(i,i+2)$, where $i\equiv 0,1 \pmod 4$. Set $t_0=\pi(m-1)$. By the Dirichlet principle, we can find an index $i$ such that $(n+i)(n+i+2)$ is divisible by at most
$$N= \left\lfloor \frac {t-t_0} {\lfloor \frac m 4 \rfloor+ \lfloor \frac {m+1} 4 \rfloor} \right\rfloor$$ of the primes from the set
\begin{equation}
\label{set}
\{p_{t_0+1}, \ldots, p_t\}.
\end{equation}
Now one writes $X=n+i+1$ and $Y=(n+i)(n+i+2)$, and gets $$X^2-DY^2=1,$$
where $D$ is squarefree and divisible only by some of the first $t_0$ primes and at most $N$ of the primes from the set (\ref{set}). Then one proceeds exactly as Lehmer, with the difference that in this approach $P(X)\leq p_t$ also has to hold. This lowers the number of Pell equations one needs to consider from $2^t-1$ to
$$M=-1+2^{t_0}\sum_{j=0}^{N}\binom{t-t_0}{j},$$
and equally important, reduces the size of the Pell equations.

The bottleneck of both the algortihms of Lehmer and Bauer and Bennett is solving the Pell equation. The main difficulty in solving the Pell equation is the size of the solutions, as it grows exponentially with respect to the size of the coefficient $d$ form (\ref{pell}). This means that just \emph{writing down} the solution in standard representation takes exponential time. In \cite{ln}, Luca and the author used \emph{compact representations} of the solutions to the Pell equations to overcome this difficulty and managed to extend Lehmer's results (see also \cite{fn2} for another application of this approach). Note that the algorithm described in \cite{ln} is still the best up to date if one wants to find 2 or 3 consecutive smooth integers.

 A compact representation of an algebraic number $\beta \in \mathbb Q (\sqrt d)$ is a representation of $\beta$ of the form
\begin{equation}
\beta=\prod_{j=1}^k \alpha_j^{2^{j}},
\label{cr}
\end{equation}
where $\alpha_j=(a_j+b_j\sqrt d)/c_j,$ where $a_j,b_j,c_j, k\in \mathbb Z$, and there exist polynomial (in $\log d$) upper bounds for $k$, $a_j,\ b_j$ and $c_j$. A detailed description of compact representations and their use can be found in \cite{jw}. Using compact representations, the Pell equation is solved in two steps. First, the regulator of the appropriate quadratic field is computed. Once the regulator is known, one computes a compact representation, in polynomial time. By far the harder part of this process is the computation of the regulator. The fastest unconditional algorithms are still exponential with respect to $d$, while the fastest known algorithm, Buchmann's subexponential algorithm (see \cite{bc}), depends on the Generalized Riemann Hypothesis.

For our purposes, the only algorithm fast enough is Buchmann's algorithm. We will perform a check, using continued fractions, that will for each case tell us either that the output of Buchmann's algorithm is unconditionally correct or that we can disregard this case. This implies that, although we cannot unconditionally solve the Pell equation in subexponential time, we can determine unconditionally whether it has a smooth solution (and find any smooth solutions) in subexponential time.

We do this by the following method, first used in \cite{ln}: suppose we want to find the solutions of $x^2-dy^2=1$, where $x$ and $y$ are $n$-smooth ($d$ is a fixed $n$-smooth integer).
We first compute the regulator $R_d$ of $\mathbb Q(\sqrt{d})$ using Buchmann's subexponential algorithm.
As mentioned, the correctness of the output $R_d$ is dependent on the GRH, but what is unconditionally true is that the output is equal to $mR_d$, where $m$ is some positive integer. From the obtained value we construct, in polynomial time, a compact representation of the $m$-th power of the fundamental unit, and from it a compact representation of the $m$-th solution $x_m+y_m\sqrt d$ of our Pell equation. We define $z$ to be the $n$-smooth part of $y_m$, which is computed using the algorithm for performing modular arithmetic on compact representations from \cite{fn}.

Let $x_1+y_1\sqrt d$ be the fundamental solution of our Pell equation and $x_k+y_k\sqrt d=(x_1+y_1\sqrt d)^k$. If the GRH is false and  Buchmann's algorithm does not return the exact value of the regulator this will mean that it might happen that we will miss some smooth solutions of our equation.  Suppose this happens. This means that $m>1$ and $y_k$ is $n$-smooth for some $k$ that is not a multiple of $m$. But we have $y_1|y_k$, so $y_1$ is also $n$-smooth, implying $y_1|z$ (because $y_1|y_m$). This obviously implies $y_1 \leq z$. So if something has gone wrong, we will find the solution $y_1$ that is smaller or equal to $z$ and strictly smaller than $y_m$. We can use continued fractions to check whether such a $y_1$ exists. If we do not find such a $y_1$ (and indeed we never do, as it's existence would disprove the GRH), it means either $m=1$ or $m>1$ but we are not missing a smooth solution. In any case, this removes the dependence of our results on the Generalized Riemann Hypothesis. Note that although $y_m$ does become very large, $z$ in practice does not. This means that it is usually, for our purposes, enough to compute the first 10 or so  convergents of the expansion of $\sqrt d$ to a continued fraction, to prove that the hypothetical $y_1$ does not exist.

Using compact representation, in combination with this check, cuts down the space needed from exponential to polynomial with respect to $d$ and the time needed to solve the Pell equation from exponential to subexponential.

In essence, in this paper we combine the clever approach of \cite{bb} and the powerful methods for finding smooth solutions of the Pell equations from \cite{ln}.

It is worth noting the somewhat counter--intuitive implications of the mentioned algorithms: it is easier to solve the Pell equation (using compact representations) than it is to write down the solution in decimal notation, and it is easier to unconditionally find all smooth solutions of the Pell equation or determine that there are none (this takes subexponential time) than it is to unconditionally solve the Pell equation (this takes exponential time).

\section{Results}
We run our algorithm for two pairs of parameters $(m,t)$, these pairs being $(47,14),$ and $(56,16)$. Note thet $p_{47}=211,$ $p_{56}=263$ and $P(\Pi_{318, 13})=163,$ $P(\Pi_{1330, 15})=223$, implying
\begin{equation}
f(k)\geq 14 \text{ for }178 \leq k\leq 223\text{ and }f(k)\geq 16\text{ for }224 \leq k\leq 268
\end{equation}
By proving that there is no integer $n>k$ for $178 \leq k\leq 222$ satisfying $P(\Pi_{n, 14})\leq k$ implies $f(k)=14$ for $178 \leq k\leq 222$. In the same way one proves that $f(k)= 16$ for $223 \leq k\leq 268$.

The case $(m,t)=(47,14)$ makes us solve
$$M_1=-1+2^{6}\sum_{j=0}^{6}\binom{41}{j}=342948991$$
Pell equations, while the case $(m,t)=(56,16)$ makes us solve
$$M_2=-1+2^{6}\sum_{j=0}^{6}\binom{50}{j}=1168680703$$
Pell equations.

We obtain that there do not exist such integers $n$, thus proving Theorem 1.\\

\textbf{Acknowledgements.}\\
These materials are based on work financed by the National Foundation for Science, Higher Education and Technological Development of the Republic of Croatia.

\small{MATHEMATISCH INSTITUUT, P.O. BOX 9512, 2300 RA LEIDEN, THE NETHERLANDS}\\

\emph{E-mail address:} fnajman@math.leidenuniv.nl\\

AND\\

\small{DEPARTMENT OF MATHEMATICS, UNIVERISTY OF ZAGREB, BIJENI\v CKA CESTA 30, 10000 ZAGREB, CROATIA}\\

\emph{E-mail address:} fnajman@math.hr

\end{document}